\documentclass[11pt,leqno]{article}
\usepackage[utf8]{inputenc}

\usepackage{amsmath}
\usepackage{amsthm}
\usepackage[shortlabels]{enumitem} \setlist{wide, labelindent=0pt}
\usepackage{mathtools}

\usepackage[colorlinks=false]{hyperref}

\usepackage{color}

\usepackage[charter]{mathdesign}

\makeatletter
\g@addto@macro\bfseries{\boldmath}
\makeatother

\usepackage{xy}
\xyoption{all}

\numberwithin{equation}{section}

\theoremstyle{plain}
    \newtheorem{theorem}[equation]{Theorem}
    \newtheorem{lemma}[equation]{Lemma}
    \newtheorem{corollary}[equation]{Corollary}
    \newtheorem{proposition}[equation]{Proposition}

    \newtheorem*{theorem*}{Theorem}
    \newtheorem*{proposition*}{Proposition}
    \newtheorem*{corollary*}{Corollary}
    \newtheorem*{lemma*}{Lemma}
    \newtheorem*{conjecture*}{Conjecture}
    \newtheorem{definition-theorem}[equation]{Definition/Theorem}
    \newtheorem{definition-lemma}[equation]{Definition/Lemma}
\theoremstyle{definition}
    \newtheorem{definition}[equation]{Definition}

    \newtheorem{remark}[equation]{Remark}
    
    \newtheorem{remarks}[equation]{Remarks}

   	\renewcommand{\phi}{\varphi}
	\let\epsilon\varepsilon

    \newcommand{\Compact}{\operatorname{K}}
    \newcommand{\Adjointable}{\operatorname{L}}
    \newcommand{\Multiplier}{\operatorname{M}}
    \newcommand{\cb}{\mathrm{cb}}
    
    \newcommand{\alg}{\mathrm{alg}}
    \newcommand{\category}{\mathsf}
   \newcommand{\functor}{\mathcal}

\newcommand{\into}{\hookrightarrow}

\newcommand{\id}{\mathrm{id}}

  \renewcommand{\prod}{\bigsqcap}

    \DeclareMathOperator{\End}{End}

    \DeclareMathOperator{\lspan}{span}

 \newcommand{\h}{\mathrm{h}}

	\newcommand{\bra}[1]{{\langle{#1}\vert}}
	\newcommand{\ket}[1]{{\vert {#1}\rangle}}

\newcommand{\Frob}{\operatorname{Frob}}
\newcommand{\LAdj}{\operatorname{LAdj}}

\title{Frobenius $C^*$-algebras and local adjunctions of $C^*$-correspondences
 
 
}
\author{Tyrone Crisp\thanks{Department of Mathematics \& Statistics, University of Maine. \href{mailto:tyrone.crisp@maine.edu}{\texttt{tyrone.crisp@maine.edu}}}}

\date{August, 2021}

\begin{document}

\maketitle

\begin{abstract}
Generalising the well-known correspondence between two-sided adjunctions and Frobenius algebras, we establish a one-to-one correspondence between local adjunctions of $C^*$-correspondences, as defined and studied in prior  work with P.~Clare and N.~Higson; and Frobenius $C^*$-algebras, a natural $C^*$-algebraic adaptation of the standard notion of Frobenius algebras that we introduce here. 
\end{abstract}

\section{Introduction}

In the theory of $C^*$-algebras one frequently encounters algebras without a multiplicative identity (often called non-unital algebras): $C^*$-algebras of compact operators on infinite-dimensional Hilbert spaces; of functions on non-compact topological spaces; and of non-discrete groups are some prominent examples. There are, however, useful proxies for multiplicative identities---multiplier algebras  and approximate identities---which often make it possible to adapt constructions from other parts of algebra, where multiplicative identities are assumed, to the non-unital $C^*$-algebra setting. Notable examples include  Gelfand duality, where the morphisms on the $C^*$-algebra side corresponding to continuous maps of non-compact spaces are defined in terms of the multiplier algebra (see, e.g., the discussion in \cite[Section 1]{anHuef-Raeburn-Williams}); and the theory of multiplier Hopf algebras and locally compact quantum groups, in which the comultiplication maps take values in a multiplier algebra \cite{Kustermans-Vaes}. Similar techniques are available, more generally, in \emph{$C^*$-categories}; see for instance \cite{Kandelaki,Vasselli,Antoun-Voigt}.

This paper is concerned with functors between $C^*$-categories of Hilbert $C^*$-modules. For each $C^*$-algebra $A$, consider the category $\Adjointable_A$ of right Hilbert $C^*$-modules over $A$, with adjointable maps as morphisms. (See Section \ref{sec:background} for background and references on Hilbert $C^*$-modules.) There are several reasons why the most useful way to view $\Adjointable_A$ is not as a category (or a $C^*$-category) on its own, but as the \emph{multiplier category} of the non-unital $C^*$-category $\Compact_A$ of  Hilbert $C^*$-modules over $A$ and \emph{compact} operators, where `compact' has the usual $C^*$-module meaning: $\Compact_A(X,Y)=\overline{\lspan}\{\ket{y}\bra{x} \in \Adjointable_A(X,Y)\, |\, x\in X,\ y\in Y\}$. For one thing, viewing $\Adjointable_A$ as the multiplier category of $\Compact_A$ accommodates the many applications in which the distinction between compact and adjointable operators is important: e.g., in operator $K$-theory; in the study of Cuntz-Pimsner algebras and related algebras; and in the study of Jones index theory and adjunctions between $C^*$-correspondences undertaken in \cite{KPW} and \cite{CCH2}, about which we shall say more below. 

Another point in favour of the multiplier-category point of view emerges when one considers functors $\Adjointable_A \to \Adjointable_B$ between $C^*$-module categories. In the general theory of multiplier $C^*$-categories the most natural functors to consider are the (multiplier extensions of) nondegenerate $*$-functors; see \cite[p.1687]{Antoun-Voigt} for the general definition. In the case of $\Compact_A$ and $\Adjointable_A$ these functors turn out to be precisely the ones that are of most interest in $C^*$-module theory: tensor products with $C^*$-correspondences. In particular, the notion of equivalence coming from the theory of multiplier $C^*$-categories corresponds, in the case of  the categories $\Compact_A$ and $\Adjointable_A$, to the usual notion of (strong) Morita equivalence of $C^*$-algebras. These results are due to Blecher \cite[Theorems 5.4 \& 5.5]{Blecher-newapproach}; cf. \cite[Proposition 2.7]{Antoun-Voigt} for a related result on a slightly different category.

If the framework of multiplier $C^*$-categories is better suited to the categorical study of Hilbert $C^*$-modules than is the standard  theory of categories, then we should expect that standard categorical notions might need to be adapted to this new framework in order to make them most useful in applications. This paper is based on a modification of the notion of adjoint functors between categories of $C^*$-modules that was introduced in \cite{CCH2}. Let $\functor{F}:\Adjointable_A\to \Adjointable_B$ and $\functor{E}:\Adjointable_B\to \Adjointable_A$ be nondegenerate $*$-functors (i.e., functors of tensor product with $C^*$-correspondences). An \emph{adjunction}, in the standard sense of category theory, between $\functor{F}$ and $\functor{E}$ is a family of natural isomorphisms
\begin{equation}\label{eq:adj-intro}
\Adjointable_B( X, \functor{F}Y) \cong \Adjointable_A(\functor{E}X,Y),
\end{equation}
one for each Hilbert $B$-module $X$ and each Hilbert $A$-module $Y$. Note that such adjunctions are always \emph{two-sided} (or \emph{ambidextrous}), since we can apply the $*$ operation on each side of \eqref{eq:adj-intro} to obtain natural isomorphisms $\Adjointable_B(\functor{F}Y,X)\cong \Adjointable_A(\functor{E}X,Y)$. (These are examples of \emph{dagger adjunctions}, cf. \cite{Heunen-Karvonen}.) Kajiwara, Pinzari, and Watatani studied adjunctions of this sort in \cite[Section 4]{KPW}, and showed that adjunctions can be characterised by a finite-index condition on $C^*$-correspondences. This finiteness condition is quite restrictive, and in order to apply adjoint functor techniques to a wider range of examples a weaker notion of adjunction  was introduced  in \cite{CCH2}: a \emph{local adjunction} between  $\functor{F}:\Adjointable_A\to \Adjointable_B$ and $\functor{E}:\Adjointable_B \to \Adjointable_A$ is a family of natural isomorphisms
\begin{equation}\label{eq:ladj-intro}
\Compact_B( X, \functor{F}Y) \cong \Compact_A(\functor{E}X,Y)
\end{equation}
of morphism spaces in the categories $\Compact_B$ and $\Compact_A$, that need not extend to the multiplier categories (that is, to an adjunction as in \eqref{eq:adj-intro}). Covering maps furnish a good illustration of the difference between adjunctions and local adjunctions: finite covers of compact Hausdorff topological spaces give rise to adjunctions, while finite \emph{branched} covers give rise to local adjunctions (cf. \cite{Pavlov-Troitsky}). Local adjunctions, like adjunctions, are automatically two-sided: from \eqref{eq:ladj-intro} we also obtain natural isomorphisms  $\Compact_B( \functor{F}Y,X) \cong \Compact_A(Y,\functor{E}X)$ by applying $*$ on both sides.

Two-sided adjunctions are known to correspond, in the very general setting of adjunctions in arbitrary $2$-categories,  to \emph{Frobenius algebras} in monoidal categories; see \cite{Morita-adjointpairs}, \cite{Muger-Frobenius}, \cite{Street-Frobenius}, and \cite{Lauda} for different aspects of this correspondence, and see \cite{Heunen-Karvonen} for analogous results for dagger adjunctions. The adjunctions studied in \cite{KPW} fit into the latter framework, and the corresponding Frobenius algebras have recently been investigated, in the case of  bimodules over a commutative $C^*$-algebra, in \cite{Heunen-Reyes}.

In this paper we show that  the connection between two-sided adjunctions and Frobenius algebras extends to the  setting of local adjunctions between $C^*$-correspondences. We introduce a notion of \emph{Frobenius $C^*$-algebra} that is not an instance of the standard algebraic definition, but  is a natural adaptation of the standard definition to the $C^*$-algebraic setting, in which conditions related to multiplicative identities are weakened using multiplier algebras and approximate identities, and the algebraic tensor product is replaced by an appropriate operator-theoretic tensor product. We  prove that isomorphism classes of Frobenius $C^*$-algebras correspond to isomorphism classes of local adjunctions, in a way that extends the known correspondence between Frobenius algebras and adjunctions. These results are close analogues of well-known theorems of Morita in the algebraic setting \cite{Morita-adjointpairs}, but the proofs are substantially different: the main step in our argument is to establish an equivalence between two natural families of norms on the spaces of square matrices over a Frobenius $C^*$-algebra.

Our results are not obviously a special case of those of \cite{Lauda} or \cite{Heunen-Karvonen}, but it would be very interesting to realise them as such: that is,  to find a $2$-category $\category{C}$ (or an appropriate $C^*$-algebraic analogue of a $2$-category)  such that the $1$-equivalences in $\category{C}$ are precisely the Morita equivalences of $C^*$-algebras, and the two-sided (or dagger) adjunctions in $\category{C}$ are precisely the local adjunctions of nondegenerate $*$-functors, yielding the correspondence established here between local adjunctions and Frobenius $C^*$-algebras as an instance of the general correspondence between two-sided adjunctions and Frobenius algebras. Such a $2$-category would be a strong contender for the `correct' categorical framework in which to study non-unital $C^*$-algebras and their modules, and would potentially provide an interesting new context in which to study categorical quantum mechanics, complementing \cite{Abramsky-Heunen} and \cite{Heunen-Reyes}, for instance. The challenge in this program is to identify the appropriate class of $2$-morphisms for the hypothetical $2$-category $\category{C}$, since the would-be unit/counit maps in a local adjunction are not adjointable (cf.~\cite[Section 3.5]{CCH2}).  In any case, we contend that the correspondence between local adjunctions and Frobenius-type objects that is established in this paper furnishes positive evidence for the existence of such a $2$-category. We leave the further pursuit of $\category{C}$, and the investigation of Frobenius-type objects in other $C^*$-categories, to future work.

\paragraph*{Plan of the paper.} In Section \ref{sec:def} we recall the standard definition of Frobenius algebras and present our $C^*$-algebraic adaptation of the definition. In Section \ref{sec:ladj-to-frob} we review some basic facts about local adjunctions of $C^*$-correspondences, and then prove that every local adjunction gives rise to a Frobenius $C^*$-algebra. In Section \ref{sec:frob-to-ladj} we reverse the construction, obtaining a local adjunction from each Frobenius $C^*$-algebra, and in Section \ref{sec:isos} we show that our constructions give a perfect correspondence between isomorphism classes of Frobenius $C^*$-algebras and of local adjunctions. Expository material and references on those parts of  $C^*$-algebra theory necessary for reading the paper are interspersed throughout.

\section{Frobenius algebras and $C^*$-algebras}\label{sec:def}

In the algebraic study of Frobenius algebras (by which we mean what are often called Frobenius extensions, as in  \cite{Kadison-Frobenius} for example), all rings  are  assumed to have multiplicative identities, and this assumption is quite important. This is particularly evident in the following formulation of the  definition:

\begin{definition}\label{def:Frobenius-algebra}
Let $A$ be a ring with multiplicative identity. A \emph{Frobenius algebra over $A$} is a ring  $C$ with multiplicative identity, equipped with an identity-preserving ring homomorphism $\eta:A\to C$ (which we use to regard $C$ as an $A$-bimodule), and an $A$-bimodule map $\epsilon:C\to A$ such that the ring $C\otimes^{\alg}_A C$ with multiplication defined by $(c_1\otimes c_2)\cdot (c_3\otimes c_4)\coloneqq c_1\epsilon(c_2 c_3)\otimes c_4$ has a multiplicative identity.
\end{definition}

Here $\otimes_A^{\alg}$ denotes the standard algebraic tensor product of $A$-modules. There are many other common definitions of Frobenius algebras (extensions), all equivalent to Definition \ref{def:Frobenius-algebra}. For instance, an element $\sum_{i=1}^n x_i\otimes y_i \in C\otimes^{\alg}_A C$ is a multiplicative identity if and only if $\{(x_i,y_i)\, |\, i=1,\ldots,n\}$ is a quasi-basis  for $\epsilon$ in the sense of \cite[Definition 1.2.2]{Watatani}.

Applying standard $C^*$-algebraic techniques for dealing with non-unital algebras, and replacing the algebraic tensor product by a tensor product suited to modules over $C^*$-algebras, we arrive at the following adaptation of Definition \ref{def:Frobenius-algebra}. Some of the terminology and notation will be explained below.

\begin{definition}\label{def:Frobenius-Cstar}
Let $A$ be a $C^*$-algebra. A \emph{Frobenius $C^*$-algebra over $A$} is a $C^*$-algebra $C$ equipped with a  $*$-homomorphism $\eta:A\to \Multiplier(C)$ (which we use to regard $C$ as an $A$-bimodule) and a positive $A$-bimodule map $\epsilon:C\to A$ such that
\begin{enumerate}[\rm(1)]
\item the algebra $C\otimes^{\h}_A C$ with multiplication $(c_1\otimes c_2)\cdot (c_3\otimes c_4)\coloneqq c_1\epsilon(c_2 c_3)\otimes c_4$ has a bounded approximate identity; and
\item $C=\left\{\eta(a)c \, \left|\, a\in \overline{\epsilon(C)},\ c\in C \right. \right\}$.
\end{enumerate}
\end{definition}

\begin{remarks}
\begin{enumerate}[\rm(1)]
\item $\Multiplier(C)$ denotes the \emph{multiplier algebra} of $C$: that is, the $C^*$-algebra of linear maps $t:C\to C$ for which there exists a linear map $t^*:C\to C$ satisfying $t(c_1)^* c_2=c_1^* t^*(c_2)$ for all $c_1,c_2\in C$. There is a canonical embedding of $C^*$-algebras $C\into \Multiplier(C)$, sending $c\in C$ to the operator $c'\mapsto cc'$. See \cite[II.7.3]{Blackadar}.
\item $\otimes^{\h}_A$ denotes the \emph{Haagerup tensor product} of operator $A$-modules, a Banach-space completion of the algebraic tensor product $\otimes^{\alg}_A$. The Haagerup tensor product is in fact not just a Banach space, but an \emph{operator space}, meaning in particular that for each $n\geq 1$ the space $M_n(C\otimes^{\h}_A C)$ of $n\times n$ matrices over $C\otimes^{\h}_A C$ is equipped with a canonical Banach-space norm.  See \cite[Section 3.4]{BLM} for the definition of the Haagerup tensor product, the details of which will not be essential for reading this paper. The appearance of the Haagerup tensor product is natural here, since this is precisely the tensor product that linearises multiplication in $C^*$-algebras (cf. \cite[Theorem 1.5.7]{BLM}); and it is also the tensor product that links $C^*$-correspondences with functors on $C^*$-module categories, as we shall recall in Section \ref{sec:background}.\item A map $\epsilon:C\to A$ is \emph{positive} if it sends positive elements to positive elements: that is, if for all $c\in C$ there is an $a\in A$ such that $\epsilon(c^*c)=a^*a$.
\item If $\epsilon:C\to A$ is a positive $A$-bimodule map, then it is automatically \emph{completely positive}: that is, for each $n\geq 1$ the map on $n\times n$ matrices $\epsilon:M_n(C)\to M_n(A)$ defined by applying $\epsilon$ to each matrix entry is also positive (see e.g. \cite[II.6.10.2]{Blackadar}). This implies in particular that $\epsilon$ is \emph{completely bounded}:  the quantity $\|\epsilon\|_{\cb}\coloneqq \sup_{n\geq 1} \| \epsilon:M_n(A)\to M_n(C)\|_{\text{operator}} $ is finite (see e.g. \cite[1.3.3]{BLM}). The complete boundedness of $\epsilon$, and of the multiplication in the $C^*$-algebra $\Multiplier(C)$, ensures that the given formula for the multiplication on the algebraic tensor product $C\otimes^{\alg}_A C$ extends to the Haagerup tensor product $C\otimes^{\h}_A C$: see \cite[3.4.5]{BLM}.
\item A \emph{bounded approximate identity} for a Banach algebra $B$  is a bounded net $\{x_\lambda\in B\, |\, \lambda\in \Lambda\}$ ($\Lambda$ being some directed set) such that  for all $y\in B$ one has $y = \lim_\lambda x_\lambda y = \lim_\lambda yx_\lambda$. Every $C^*$-algebra has a bounded approximate identity: see e.g. \cite[II.4.1]{Blackadar}.
\item If the $C^*$-algebra $C$ has a multiplicative identity then condition (2) in the definition above is implied by condition (1): see Remark \ref{remark:C-unital}.
\item The Cohen-Hewitt factorisation theorem (cf. \cite[II.5.3.7]{Blackadar}) ensures that the right-hand side in Condition (2) is equal to $\overline{\lspan}\{ \eta\circ \epsilon(c_1)c_2\, |\, c_1,c_2\in C\}$. Since we shall tacitly invoke the factorisation theorem at several points in this paper, let us recall what it says: if $A$ is a $C^*$-algebra, and if $X$ is a Banach space and a left $A$-module for which there is a positive constant $c$ satisfying $\|ax\|\leq c \|a\| \|x\|$ for every $a\in A$ and $x\in X$, then $\{ax\in X\ |\ a\in A,\ x\in X\}$ is a closed $A$-submodule of $X$. 
\item Even if the $C^*$-algebras $A$ and $C$ have multiplicative identities, and the homomorphism $\eta:A\to \Multiplier(C)=C$ satisfies $\eta(1_A)=1_C$, condition (1) is still strictly weaker than the requirement that $C\otimes^{\h}_A C$ have a multiplicative identity. For example, take $C=C[-1,1]$ (continuous, complex-valued functions on the real interval $[-1,1]$). For each $c\in C$ let $c^{-}\in C$ be the function $x\mapsto c(-x)$, and let $A\subseteq C$ be the subalgebra $\{a\in C\, |\, a=a^{-}\}$ of even functions. Consider the inclusion map $\eta:A\to C$, and the averaging map $\epsilon:C\to A$ given by $\epsilon(c)\coloneqq \frac{1}{2}(c+c^{-})$. Then the map 
\[
c_1\otimes c_2\mapsto \frac{1}{2}\begin{bmatrix} c_1c_2 & c_1 c_2^{-} \\ c_1^{-}c_2 & c_1^{-}c_2^{-} \end{bmatrix}
\]
gives a Banach-algebra isomorphism from $C\otimes^{\h}_A C$ to the $C^*$-algebra
\[
\left\{ \left. \begin{bmatrix} c_1 & c_2 \\ c_2^{-} & c_1^{-}\end{bmatrix}\in M_2(C)\, \right|\, c_1,c_2\in C,\ c_1(0)=c_2(0)\right\}.
\]
The latter algebra does not possess a multiplicative identity; but, like all $C^*$-algebras, it does possess a bounded approximate identity. (This is a standard example in the study of different notions of finite index for inclusions of $C^*$-algebras: see for instance \cite[Proposition 2.8.2]{Watatani}, \cite[p.88]{Frank-Kirchberg}, \cite[Example 2.32]{KPW}.)
\item An alternative definition of Frobenius algebras in the algebraic setting involves a comultiplication map $\delta:C\to C\otimes^{\alg}_A C$ \cite{Abrams-Frobenius}. A similar formulation is possible for Frobenius $C^*$-algebras, with the comultiplication being a map $\delta: C\to \Multiplier(C\otimes^{\h}_A C)$; this is reminiscent of the way comultiplication is treated for locally compact quantum groups in \cite{Kustermans-Vaes}. Some care is necessary, however, in formulating Definition \ref{def:Frobenius-Cstar} in this way, since one must first define  $\Multiplier(C\otimes^{\h}_A C)$. If we insist that $C\otimes^{\h}_A C$ be a $C^*$-algebra then the meaning is clear, but that assumption already implies that $C$ is a Frobenius $C^*$-algebra (see remark (5)), and the existence of a comultiplication $C\to \Multiplier(C\otimes^{\h}_A C)$ is then guaranteed: take a bounded approximate identity $\{x_\lambda\}$ for $C\otimes^{\h}_A C$ and then define $\delta(c)\coloneqq \operatorname{strict-lim}_\lambda c x_\lambda$, the limit in the strict topology of the multiplier algebra.  One can perhaps avoid this redundancy by referring to the more general notion of multipliers of operator spaces (as in \cite[Section 4.5]{BLM}), but it is not clear to us how to neatly express in that framework the condition that the map $\epsilon$ should act as a counit for the comultiplication. 
\end{enumerate}
\end{remarks}

\section{Frobenius $C^*$-algebras from local adjunctions}\label{sec:ladj-to-frob}

We begin this section by recalling the algebraic model upon which our $C^*$-algebraic result is based. Let $A$ and $B$ be rings with $1$, and let ${}_A F_B$ be an $A$-$B$-bimodule such that the functor $X\mapsto X\otimes^{\alg}_A F$, from the category of right $A$-modules to the category of right $B$-modules, has a two-sided adjoint. Then the ring $C\coloneqq \End_B(F)$ is a Frobenius algebra over $A$ in the sense of Definition \ref{def:Frobenius-algebra}: see \cite[Theorem 9.3]{Morita-adjointpairs}. The main result of this section, Proposition \ref{prop:ladj-to-frob}, shows that an analogous construction produces Frobenius $C^*$-algebras from local adjunctions of $C^*$-correspondences.

\subsection{Background on Hilbert $C^*$-modules and local adjunctions}\label{sec:background}

Let us first establish some basic notions and notation. Let $A$ be a $C^*$-algebra (\cite{Blackadar} is a convenient general reference). A \emph{Hilbert $A$-module} is, for us, a right $A$-module $X$, equipped with an $A$-valued inner product $\langle\, |\, \rangle : X\times X\to A$ satisfying a list of axioms analogous to those of a complex Hilbert space (which is the same thing as a Hilbert $\mathbb{C}$-module). In particular, the norm $x\mapsto \| \langle x\, |\, x\rangle \| ^{1/2}$ makes $X$ into a Banach space. For details see, e.g., \cite{Lance}, \cite[II.7.1]{Blackadar}, or \cite[Chapter 8]{BLM}. 

If $X$ is a Hilbert $A$-module then the space $M_n(X)$  of $n\times n$ matrices over $X$ is in a natural way a Hilbert module over the $C^*$-algebra $M_n(A)$, and $M_n(X)$ thus comes equipped with a canonical norm. Using these norms, the notions of completely bounded, completely contractive, and completely isometric maps, and also the Haagerup tensor product, make sense for Hilbert $C^*$-modules. See \cite[Section 8.2]{BLM} for details.

An $A$-linear map $t: X\to Y$ of Hilbert $A$-modules is called \emph{adjointable} if there is a map $t^*:Y\to X$ satisfying $\langle t(x)\, |\, y\rangle = \langle x \, |\, t^*(y)\rangle$ for all $x\in X$ and $y\in Y$. Every adjointable map is completely bounded. We denote by $\Adjointable_A$ the category whose objects are Hilbert $A$-modules, and whose morphisms are the adjointable maps. This is a \emph{$C^*$-category} as explained, for instance, in \cite{Antoun-Voigt}. In particular, for each Hilbert $A$-module $X$ the endomorphism algebra $\Adjointable_A(X) =\Adjointable_A(X,X)$ is a $C^*$-algebra.

An important subclass of the adjointable maps are the \emph{compact} ones. For each $x\in X$ we let $\bra{x}: A\to X$ be the map $a\mapsto xa$. This is an adjointable map, whose adjoint $\bra{x}\coloneqq \ket{x}^*$ is given by $x'\mapsto \langle x\,  |\, x'\rangle$. Now for each pair of Hilbert $A$-modules $X$ and $Y$ we define
\[
\Compact_A(X,Y) \coloneqq \overline{\lspan}\big\{ \ket{y}\bra{x}\in \Adjointable_A(X,Y)\, \big|\, x\in X,\ y\in Y \big\},
\]
called  the space of compact operators from $X$ to $Y$. (Note that these operators are typically not compact in the sense of mapping bounded sets to precompact sets.) The collection of compact operators forms a closed two-sided $*$-ideal in the category $\Adjointable_A $, and in particular $\Compact_A(X)=\Compact_A(X,X)$ is a closed, two-sided ideal in the $C^*$-algebra $\Adjointable_A(X)$, for every Hilbert $A$-module $X$. (The category $\Adjointable_A$ can be identified with the \emph{multiplier category} of the non-unital $C^*$-category $\Compact_A$. The  theory of nonunital $C^*$-categories and their multiplier categories is an important source of motivation for this paper, but we shall not need the details here; for that see, e.g., \cite{Antoun-Voigt}.) 

The next theorem recalls some of the important properties of compact operators that we shall use:

\begin{theorem}\label{thm:compacts}
Let $X$ be a Hilbert $A$-module.
\begin{enumerate}[\rm(1)]
\item The map sending $t\in \Adjointable_A(X)$ to the map $\Compact_A(X)\xrightarrow{k\mapsto tk} \Compact_A(X)$ gives an isomorphism of $C^*$-algebras $\Adjointable_A(X)\cong \Multiplier(\Compact_A(X))$.
\item When $A$ is regarded as a Hilbert $A$-module, with the obvious module structure and with the inner product $\langle a_1\,  |\, a_2\rangle \coloneqq a_1^* a_2$, one has $A\cong \Compact_A(A)$ via the map sending $a\in A$ to the operator $\ket{a}$.
\item The map $\Compact_A(A,X)\otimes^{\h}_A \Compact_A(X,A)\xrightarrow{k_1\otimes k_2\mapsto k_1 k_2} \Compact_A(X)$ is a completely isometric isomorphism.
\item $X = \{k(x)\, |\, k\in \Compact_A(X),\ x\in X\}$.
\end{enumerate}
\end{theorem}

\begin{proof}
For (1) see \cite[II.7.3.10]{Blackadar}. For (2) see  \cite[II.7.2.6(i)]{Blackadar}. For (3) see \cite[8.1.11 \& 8.2.15]{BLM}, and for part (4) see \cite[Proposition 1.3]{Lance}.
\end{proof}

We now turn to functors between categories of Hilbert $C^*$-modules. For the reasons explained in the introduction (cf. Remark \ref{remark:functors} below) we shall restrict our attention to the functors of tensor product with $C^*$-correspondences. For $C^*$-algebras $A$ and $B$, a \emph{$C^*$-correspondence} from $A$ to $B$ is a Hilbert $B$-module $F$ equipped with a homomorphism of $C^*$-algebras $\eta_F:A\to \Adjointable_B(F)$, that is nondegenerate in the sense that $\Compact_B(F)=\{\eta_F(a)k\ |\ a\in A,\ k\in \Compact_B(F)\}$ (or, equivalently, $F=\{\eta_F(a)(f)\ |\ a\in A,\ f\in F\}$). We shall often drop  $\eta_F$ from the notation and just write $ak$ and $af$ for $\eta_F(a) k$ and $\eta_F(a)(f)$. We shall also write `let ${}_A F_B$ be a $C^*$-correspondence' to mean that $F$ is a $C^*$-correspondence from $A$ to $B$.

If $Y$ is a Hilbert $A$-module, and if $F$ is a $C^*$-correspondence from $A$ to $B$, then the Haagerup tensor product $Y\otimes^{\h}_A F$ is a Hilbert $B$-module: the inner product is given by
\[
\langle y_1\otimes f_1\, |\, y_2\otimes f_2\rangle = \langle f_1\, |\, \eta_F(\langle y_1\, |\, y_2\rangle)f_2\rangle,
\]
and the norm induced by this inner product on $Y\otimes^{\h}_A F$ (and on its matrix spaces) coincides with the Haagerup norm; see \cite[8.2.11]{BLM}. For each $t\in \Adjointable_A(Y,Y')$ the map $t\otimes\id_F: Y\otimes^{\h}_A F\to Y'\otimes^{\h}_A F$ lies in $\Adjointable_B(Y\otimes^{\h}_A F, Y'\otimes^{\h}_A F)$, and  tensor product with $F$ determines in this way a functor $\functor{F}:\Adjointable_A \to \Adjointable_B  $.

\begin{remark}\label{remark:functors}
The functor $\functor{F}:\Adjointable_A \to \Adjointable_B  $, $Y\mapsto Y\otimes^{\h}_A F$ of tensor product with a $C^*$-correspondence ${}_A F_B$ has the following two properties:
\begin{enumerate}[\rm(1)]
\item it is a \emph{$*$-functor}, meaning that the maps $\functor{F}:\Adjointable_A(Y,Y')\to \Adjointable_B(\functor{F}Y,\functor{F}Y')$ are linear and satisfy $\functor{F}(t^*)=\functor{F}(t)^*$; and
\item it is \emph{nondegenerate}, in the sense that for each Hilbert $A$-module $Y$ we have $\Compact_B(\functor{F}Y)=\{\functor{F}(k)\circ l\ |\ k\in \Compact_A(Y),\ l\in \Compact_B(\functor{F}Y)\}$. (An equivalent formulation is that the maps $\functor{F}:\Adjointable_A(Y,Y')\to \Adjointable_B(\functor{F}Y,\functor{F}Y')$ are strongly continuous on bounded subsets.) 
\end{enumerate}
Blecher has shown that every nondegenerate $*$-functor $\functor{F}:\Adjointable_A \to \Adjointable_B$ is canonically, unitarily naturally isomorphic to the functor $Y\mapsto Y\otimes^{\h}_A \functor{F}A$; see \cite[Theorem 5.4]{Blecher-newapproach}.
\end{remark}

We now recall the definition of local adjunctions of $C^*$-correspondences from \cite{CCH2}. The definition below is the one that will be most useful in this paper; see Remarks \ref{remark:functors2}(2) for an equivalent formulation that looks more like an adjunction in the usual categorical sense.

\begin{definition}\label{def:ladj}
A \emph{local adjunction} between $C^*$-correspondences ${}_A F_B$ and ${}_B E_A$  is a completely bounded conjugate-linear isomorphism $\phi:F\to E$ satisfying $\phi(afb)=b^*\phi(f)a^*$ for all $a\in A$, $b\in B$, and $f\in F$.
\end{definition}

\begin{remarks}\label{remark:functors2}
\begin{enumerate}[\rm(1)]
\item A conjugate-linear map $\phi:F\to E$ is called completely bounded if the operator norms of the matrix extensions $\phi: M_n(F)\to M_n(E)$, $\phi(f)_{i,j}\coloneqq \phi(f_{j,i})$ (note the transposition) are uniformly bounded. This is equivalent to the condition that $f\mapsto \phi(f)^\star$ be a completely bounded linear map from $F$ to the \emph{adjoint operator space} $E^\star$; cf \cite[1.2.25]{BLM}. The map $\phi$ is a completely bounded conjugate-linear isomorphism if it has a completely bounded conjugate-linear inverse. The involution $a\mapsto a^*$ on a $C^*$-algebra is the most prominent example of a completely bounded (indeed, competely isometric) conjugate-linear isomorphism.
\item If $\phi:F\to E$ is a local adjunction in the above sense, then the maps
\[
\Phi_{X,Y}:\Compact_B(X,Y\otimes^{\h}_A F) \to \Compact_A(X\otimes^{\h}_B E, Y),\qquad \ket{y\otimes f}\bra{x}\mapsto \ket{y}\bra{x\otimes \phi(f)}
\]
are completely bounded natural isomorphisms. Conversely, each collection of natural isomorphisms of this kind comes from a local adjunction as in Definition \ref{def:ladj}; see \cite[Theorem 3.10]{CCH2}. Note that local adjunctions are automatically two-sided: the formula $\ket{x}\bra{y\otimes f}\mapsto \ket{x\otimes \phi(f)}\bra{y}$ defines an isomorphism $\Phi_{X,Y}^{\dagger} :\Compact_B(Y\otimes^{\h}_A F, X)\to \Compact_A(Y,X\otimes^{\h}_B E)$.
\item If $\phi:F\to E$ is a local adjunction then equipping $F$ with the $A$-valued inner product $\langle f_1\, |\, f_2\rangle \coloneqq \langle \phi(f_1)\, |\, \phi(f_2)\rangle$ makes $F$ into a \emph{bi-Hilbertian bimodule of finite numerical index} as defined in \cite{KPW}. Conversely, every bi-Hilbertian bimodule of finite numerical index determines a local adjunction; see \cite[Lemma 3.38]{CCH2}.
\item Several naturally occurring examples of local adjunctions are presented in \cite{CCH2}; let us recall two here. Firstly, conditional expectations of finite index in the sense of \cite{Frank-Kirchberg} (that is, those satisfying a Pimsner-Popa-type inequality) give rise to local adjunctions; by contrast, conditional expectations of finite index in the sense of \cite{Watatani} (that is, those admitting a finite quasi-basis) give rise to actual adjunctions. See \cite[Proposition 2.12]{KPW} and \cite[Section 4.3]{CCH2} for details. Secondly, if $G$ is a real reductive group and $L$ is a Levi factor of a parabolic subgroup of $G$, then there is a local adjunction between the functors of parabolic induction and restriction of Hilbert modules over the reduced group $C^*$-algebras of $G$ and $L$, and this local adjunction gives rise to an actual adjunction at the level of tempered unitary representations: see \cite{CCH1}.
\end{enumerate}
\end{remarks}

\subsection{Frobenius $C^*$-algebras from local adjunctions}

Let ${}_A F_B$ and ${}_B E_A$ be a pair of $C^*$-correspondences, and let  $\phi: F\to E$ be a local adjunction as in Definition \ref{def:ladj}. We are going to equip the $C^*$-algebra $C\coloneqq \Compact_B(F)$ with the structure of a Frobenius $C^*$-algebra over $A$.

The map $\eta:A\to C$ will be the given action of $A$ as adjointable operators on $F$:
\[
\eta_F: A \to \Adjointable_B(F)\cong \Multiplier(C).
\]
The map $\epsilon:C\to A$ will be defined by the formula
\[
\epsilon_\phi( \ket{f_1}\bra{f_2}) = \langle \phi(f_1)\, |\, \phi(f_2)\rangle.
\]
That this formula defines a positive $A$-bimodule map $\Compact_B(F)\to A$ follows from \cite[Corollary 2.11]{KPW} (which applies to the present situation thanks to \cite[Lemma 3.38]{CCH2}).

\begin{proposition}\label{prop:ladj-to-frob}
For each local adjunction $\phi: F\to E$ the $C^*$-algebra $\Compact_B(F)$, equipped with the maps $\eta_F$ and $\epsilon_\phi$, is a Frobenius $C^*$-algebra over $A$. We denote this Frobenius $C^*$-algebra by $\Frob(\phi)$. 
\end{proposition}

To prove Proposition \ref{prop:ladj-to-frob} we need to verify conditions (1) and (2) from Definition \ref{def:Frobenius-Cstar}. The nondegeneracy condition (2) is easily established, as follows. It is clear from the definition that $\overline{\epsilon_\phi(C)}= \overline{\lspan}\{ \langle e_1\, |\, e_2\rangle \in A\ |\ e_1,e_2\in E\}$.  Part (4) of Theorem \ref{thm:compacts} then implies that $E=\left\{ea\, \left|\, e\in E,\ a\in \overline{\epsilon_\phi(C)}\right. \right\}$, and pulling this equality back along the isomorphism $\phi$ yields $F=\left\{ \eta_F(a)f\, \left|\, \ a\in \overline{\epsilon_\phi(C)},\ f\in F \right.\right\}$, which immediately implies that $C=\Compact_B(F) = \left\{ \eta_F(a)c\, \left|\, a\in \overline{\epsilon_\phi(C)},\ c\in C \right. \right\}$ as required.

We are thus left to verify that the algebra $C\otimes^{\h}_A C$ has a bounded approximate identity. Since every $C^*$-algebra has a bounded approximate identity, this will follow immediately from:

\begin{lemma}\label{lem:CAC-algebra}
There is a bounded isomorphism of algebras $C\otimes^{\h}_A C\cong \Compact_A(F\otimes^{\h}_B E)$.
\end{lemma}

\begin{proof}
Consider the map $C\otimes^{\h}_A C\to \Compact_A(F\otimes^{\h}_B E)$ defined, for $f_1,f_2,f_3,f_4\in F$, by
\begin{equation}\label{eq:CAC-iso-def}
\ket{f_1}\bra{f_2}\otimes \ket{f_3}\bra{f_4} \mapsto \ket{f_1\otimes\phi(f_2)} \bra{f_4\otimes \phi(f_3)}.
\end{equation}
To see that this formula determines an isomorphism of Banach spaces (in fact, a completely bounded isomorphism of operator spaces), notice that it coincides with the composition of isomorphisms
\[
\begin{aligned}
C\otimes^{\h}_A C & = \Compact_B(F,F)\otimes^{\h}_A \Compact_B(F,F) \\
& \xrightarrow{\Phi_{F,A}^{\dagger}\otimes \Phi_{F,A}} \Compact_A(A, F\otimes^{\h}_B E)\otimes^{\h}_A \Compact_A(F\otimes^{\h}_B E,A) \\
& \xrightarrow{k_1\otimes k_2\mapsto k_1k_2} \Compact_A(F\otimes^{\h}_B E,F\otimes^{\h}_B E).
\end{aligned}
\]
(See Remarks \ref{remark:functors2} for the meaning of $\Phi_{F,A}$ and $\Phi_{F,A}^\dagger$.) A simple computation reveals that  \eqref{eq:CAC-iso-def} is an algebra map, and so it is an isomorphism of Banach algebras. 
\end{proof}

The proof of Proposition \ref{prop:ladj-to-frob} is now complete.

\begin{remark}[Adjunctions]
A local adjunction  $\phi:{}_A F_B \to {}_B E_A$ is an \emph{adjunction} if $\eta_F(A)\subseteq \Compact_B(F)$ and $\eta_E(B)\subseteq \Compact_A(E)$: that is, if $A$ acts on $F$ and $B$ acts on $E$ by compact operators. Adjunctions in this sense correspond precisely to adjunctions in the usual categorical sense \eqref{eq:adj-intro} between the tensor-product functors associated to $F$ and $E$: see \cite[Section 4]{KPW}. The Frobenius $C^*$-algebras $\Frob(\phi)$ associated to these adjunctions are precisely the \emph{dagger Frobenius monoids} (see \cite{Heunen-Karvonen}) in the dagger monoidal category of $C^*$-correspondences from $A$ to $A$, with adjointable maps as morphisms and with $\otimes^{\h}_A$ as tensor product. In the case where the $C^*$-algebra $A$ is commutative these Frobenius algebras have  been studied in \cite{Heunen-Reyes}.
\end{remark}

\section{Local adjunctions from Frobenius $C^*$-algebras}\label{sec:frob-to-ladj}

Once again we begin by recalling the algebraic model for our $C^*$-algebraic result. Suppose that $C$ is a Frobenius algebra over $A$, in the sense of Definition \ref{def:Frobenius-algebra}. Considering $C$ as an $A$-$C$-bimodule, the base-change functor $X\mapsto X\otimes^{\alg}_A C$, from right $A$-modules to right $C$-modules, has a two-sided adjoint (namely, the functor of tensor product with the bimodule ${}_C C_A$); thus we may apply the construction described at the beginning of Section \ref{sec:ladj-to-frob} to obtain another Frobenius algebra $\End_C({}_A C_C)$ over $A$, and the obvious isomorphism of rings $C\cong \End_C({}_A C_C)$ is an isomorphism of Frobenius algebras over $A$. (Cf. \cite[Theorem 5.1]{Morita-adjointpairs}.) Thus every Frobenius algebra arises from a two-sided adjunction.

We are going to prove an analogous result for Frobenius $C^*$-algebras: up to isomorphism, every Frobenius $C^*$-algebra arises from a local adjunction of $C^*$-correspondences, via the construction given in Proposition \ref{prop:ladj-to-frob}. The meaning of isomorphism for Frobenius $C^*$-algebras is the obvious one:

\begin{definition}
Let $C$ and $C'$ be two Frobenius $C^*$-algebras over $A$. An \emph{isomorphism} of Frobenius $C^*$-algebras $\rho: C \to C'$ is an isomorphism of $C^*$-algebras $\rho:C\to C'$ that is also an $A$-bimodule map, and which satisfies $\epsilon'\circ\rho=\epsilon$ (where $\epsilon:C\to A$ and $\epsilon':C'\to A$ are the Frobenius structure maps).
\end{definition}

\begin{theorem}\label{thm:frob-to-ladj}
Let $C$ be a Frobenius $C^*$-algebra over $A$. There are is a local adjunction of $C^*$-correspondences $\phi:{}_A F_B\to {}_B E_A$ such that $C \cong \Frob(\phi)$ as Frobenius $C^*$-algebras over $A$.
\end{theorem}

\subsection{Proof of Theorem \ref{thm:frob-to-ladj}}

Let $C$ be a Frobenius $C^*$-algebra over $A$, with structure maps $\eta:A\to \Multiplier(C)$ and $\epsilon:C\to A$. Consider $C={}_A C_C$ as a $C^*$-correspondence from $A$ to $C$, with the canonical right $C$-module structure and inner product $\langle c_1\, |\, c_2\rangle \coloneqq c_1^*c_2$, and with left $A$-action given by the $*$-homomorphism $\eta:A\to \Multiplier(C)=\Adjointable_C(C)$. Note that the condition (2) in Definition \ref{def:Frobenius-Cstar} ensures that $\eta$ is nondegenerate.

\begin{definition}\label{def:C-epsilon}
We construct a $C^*$-correspondence ${}_C C^{\epsilon}_A$, from $C$ to $A$, as follows. The positive $A$-bimodule map $\epsilon:C\to A$ gives an $A$-valued inner product $\langle c_1\,|\, c_2\rangle^{\epsilon} \coloneqq \epsilon(c_1^* c_2)$ on $C$, and we let $C^{\epsilon}$ denote the completion of the quotient $C/\lspan\{c\in C\, |\, \epsilon(c^*c)=0\}$ in the norm coming from this inner product. (We will later show that this competed quotient is just $C$ itself, so no quotient or completion actually takes place.) The right action of $A$ on $C$ via the homomorphism $\eta$ extends to an action on $C^{\epsilon}$, and makes $C^{\epsilon}$ a right Hilbert $A$-module. We let $q:C\to C^{\epsilon}$ denote the quotient map, and we let $C$ act on $C^{\epsilon}$ via $cq(c')\coloneqq q(cc')$. Since the multiplication action of a $C^*$-algebra on itself is nondegenerate, this turns $C^{\epsilon}$ into a $C^*$-correspondence from $C$ to $A$.
\end{definition} 

We are going to prove  that $q$ is a completely bounded isomorphism. Since $\epsilon$ is completely bounded, $q$ is likewise, and so we are left to show that $q$ is completely bounded from below. This will require some preparation.

Let $\mu:C\otimes^{\h}_A C\to C$ be the multiplication map $\mu(c_1\otimes c_2)\coloneqq c_1c_2$, and let $\epsilon^{(1)}:C\otimes^{\h}_A C\to C$ be the map $\mu\circ(\epsilon\otimes\id_C) : c_1\otimes c_2\mapsto \epsilon(c_1)c_2$. Note the following formula for the product in $C\otimes^{\h}_A C$:
\begin{equation}\label{eq:CAC-mult}
 (c_1\otimes c_2)\cdot x= c_1\otimes \epsilon^{(1)}(c_2x)\quad \text{for all}\quad c_1,c_2\in C\text{ and }x\in C\otimes^{\h}_A C.
\end{equation}

\begin{lemma}\label{lem:approx-unit-1}
Let $\{x_\lambda\, |\, \lambda\in \Lambda\}$ be a bounded approximate identity in $C\otimes^{\h}_A C$. For each $c\in C$ we have  $c = \lim_{\lambda} \epsilon^{(1)}(cx_\lambda)$.
\end{lemma}

\begin{proof}
It will suffice to prove this for elements of the form $c=\epsilon(c_1)c_2$, since such elements are dense in $C$ (condition (2) in Definition \ref{def:Frobenius-Cstar}) and the maps $c\mapsto \epsilon^{(1)}(cx_\lambda)$ are uniformly bounded. The formula \eqref{eq:CAC-mult} for the multiplication in $C\otimes^{\h}_A C$ immediately yields
\(
\epsilon^{(1)}( \epsilon(c_1)c_2 x_\lambda )=  \epsilon^{(1)}\left( (c_1\otimes c_2)\cdot x_\lambda \right),
\)
and consequently
\[
\epsilon(c_1)c_2 = \epsilon^{(1)}(c_1\otimes c_2) = \lim_\lambda \epsilon^{(1)}\left( (c_1\otimes c_2)\cdot x_\lambda\right) = \lim_{\lambda}\epsilon^{(1)}(\epsilon(c_1)c_2 x_\lambda)
\]
as required.
\end{proof}

\begin{remark}\label{remark:C-unital}
If $C$ has a multiplicative identity, then Lemma \ref{lem:approx-unit-1} can be established without appealing to condition (2) in Definition \ref{def:Frobenius-Cstar}, as follows:
\[
c = \mu(1\otimes c) = \lim_\lambda \mu\left(1\otimes \epsilon^{(1)}(cx_\lambda)\right) = \lim_\lambda \epsilon^{(1)}(cx_\lambda).
\]
Since $\epsilon^{(1)}(C\otimes^{\h}_A C) \subseteq \overline{\lspan}\{\epsilon(c_1)c_2\, |\, c_1,c_2\in C\}$, we see that condition (2) in Definition \ref{def:Frobenius-Cstar} is implied by condition (1), when $C$ has a multiplicative identity.
\end{remark}

We are now ready to prove that $q$ is completely bounded from below:

\begin{lemma}\label{lem:key-estimate} 
Let $\{x_\lambda\, |\, \lambda\in \Lambda\}$ be a bounded approximate unit for $C\otimes^{\h}_A C$, and write $s\coloneqq \sup_\lambda \|x_\lambda\|$. For each $n\geq 1$ and each $c\in M_n(C)$ we have
\[
\| c \|_{M_n(C)} \leq s \|\epsilon\|_{\cb} \|q(c)\|_{M_n(C^{\epsilon})}.
\]
\end{lemma}

\begin{proof}
Let $D=\overline{\lspan}\{\eta(a)m\, |\, a\in A,\ m\in \Multiplier(C)\}\subseteq \Multiplier(C)$. Note that $C\subseteq D$ (condition (2) in Definition \ref{def:Frobenius-Cstar}), and $\eta(A)\subseteq D$ (obviously). Equipped with the inner product $\langle d_1\, |\, d_2\rangle\coloneqq d_1^* d_2$, and the obvious $A$-$\Multiplier(C)$-bimodule structure, $D$ becomes a $C^*$-correspondence from $A$ to $\Multiplier(C)$. We shall consider the Hilbert module tensor product $C^{\epsilon}\otimes^{\h}_A D$, which is a $C^*$-correspondence from $C$ to $\Multiplier(C)$. 

The completely bounded map $q:C\to C^{\epsilon}$ and the completely isometric inclusion $\text{incl}:C\into D$ give a completely bounded map $q\otimes\text{incl}:C\otimes^{\h}_A C\to C^{\epsilon}\otimes^{\h}_A D$, and it will be useful to note that for all $c_1,c_2\in C$, all $a\in A$, and all $x\in C\otimes^{\h}_A C$ we have
\begin{equation}\label{eq:tensor-ip}
\langle q(c_1)\otimes \eta(a)\, |\, (q\otimes \text{incl})(c_2 x)\rangle = a^* \epsilon^{(1)}(c_1^*c_2x).
\end{equation}
(This formula follows easily from the definitions.) Passing to matrices, $M_n(C^{\epsilon}\otimes^{\h}_A D)$ is a Hilbert module over $M_n(\Multiplier(C))$, canonically unitarily isomorphic to $M_n(C^{\epsilon})\otimes^{\h}_{M_n(A)} M_n(D)$ (this is a straightforward variation on \cite[3.4.11]{BLM}). The inner product on $M_n(C^{\epsilon}\otimes^{\h}_A D)$ still satisfies \eqref{eq:tensor-ip}, with $c_1,c_2\in M_n(C)$, $a\in M_n(A)$, and $x\in M_n(C\otimes^{\h}_A C)$, and with $q$, $\eta$, and $\epsilon^{(1)}$ denoting the entry-wise extensions of these maps to matrices.
 
Now fix $n\geq 1$. The nondegeneracy of $C$ over $A$ implied by condition (2) in Definition \ref{def:Frobenius-Cstar} means that  every element of $M_n(C)$ can be written in the form $ca$ for some $c\in M_n(C)$ and some $a\in M_n(A)$, so it will suffice to establish the desired estimate for elements of this form. For each $\lambda\in \Lambda$ let $y_\lambda\in M_n(C\otimes^{\h}_A C)$ be the diagonal matrix with constant diagonal entries $x_\lambda$. Then $\|y_\lambda\|=\|x_\lambda\|$ for each $\lambda$, and Lemma \ref{lem:approx-unit-1} implies that $c=\lim_\lambda \epsilon^{(1)}(cy_\lambda)$ for every $c\in M_n(C)$. 

The formula \eqref{eq:tensor-ip} gives 
\[
\langle q(c)\otimes \eta(a) \, |\, (q\otimes \text{incl})(ca y_\lambda) \rangle  = a^*\epsilon^{(1)}(c^*ca y_\lambda) = \epsilon^{(1)}(a^*c^*cay_\lambda),
\]
and so 
\begin{equation}\label{eq:ip-estimate}
(ca)^*(ca) = \lim_\lambda \langle q(c)\otimes \eta(a) \, |\, (q\otimes \text{incl})(ca y_\lambda) \rangle.
\end{equation}
Now for each $\lambda\in \Lambda$ the Cauchy-Schwartz inequality  in the $C^*$-module $M_n(C^{\epsilon}\otimes^{\h}_A D)$ (see e.g. \cite[Proposition 1.1]{Lance}) gives
\begin{equation}\label{eq:CS}
\begin{aligned}
&\| \langle q(c)\otimes \eta(a) \, |\, (q\otimes \text{incl})(ca y_\lambda) \rangle\|_{M_n(C)} \\ & \leq \| q(c)\otimes \eta(a)\|_{M_n(C^{\epsilon}\otimes^{\h}_A D)}  \|q\otimes\text{incl}\|_{\cb}   \|ca\|_{M_n(C)}   \|y_\lambda\|_{M_n(C\otimes^{\h}_A C)} \\ & \leq s \|\epsilon\|_{\cb}  \|ca\|_{M_n(C)}  \|q(c)\otimes \eta(a)\|_{M_n(C^{\epsilon}\otimes^{\h}_A D)}
\end{aligned}
\end{equation}
where we have used the fact that $\|q\otimes \text{incl}\|\leq \|q\|\leq \|\epsilon\|$, and the fact that $\|y_\lambda\|=\|x_\lambda\|$. To finish the proof we note that the maps 
\[
C^{\epsilon}\otimes^{\h}_A \eta(A) \xrightarrow{\id_{C^{\epsilon}}\otimes\text{incl}} C^{\epsilon}\otimes^{\h}_A D
\quad \text{and} \quad
C^{\epsilon}\otimes^{\h}_A \eta(A)\xrightarrow{q(c)\otimes \eta(a)\mapsto q(ca)} C^{\epsilon}
\]
are completely isometric (\cite[3.6.5(2) \& 3.4.6]{BLM}), and from this it follows that 
\[
\|q(c)\otimes \eta(a)\|_{M_n(C^{\epsilon}\otimes^{\h}_A D)}= \|q(ca))\|_{M_n(C^{\epsilon})}.
\] 
Combining this last equality with \eqref{eq:ip-estimate} and \eqref{eq:CS} gives
\[
\| ca\|_{M_n(C)}^2 \leq s \|\epsilon\|_{\cb}  \|ca\|_{M_n(C)}   \|q(ca)\|_{M_n(C^{\epsilon})},
\]
and dividing throughout by $\|ca\|$ gives the estimate in the lemma.
\end{proof}

Define a map $\phi^\epsilon:C\to C^{\epsilon}$ by $\phi^\epsilon(c)\coloneqq q(c^*)$. The $*$ map on a $C^*$-algebra is a conjugate-linear complete isometry, and so Lemma \ref{lem:key-estimate} implies that the conjugate-linear map $\phi^\epsilon$ is a completely bounded isomorphism. The definition of the bimodule structures on $C$ and on $C^{\epsilon}$ make it obvious that $\phi^\epsilon(ac' c) = c^* \phi^\epsilon(c')a^*$ for all $c,c'\in C$ and all $a\in A$, and so $\phi^\epsilon$ is a local adjunction. Let us pause the proof of Theorem \ref{thm:frob-to-ladj} to make a formal definition:

\begin{definition}\label{def:LAdj}
Let $C$ be a Frobenius $C^*$-algebra over $A$, with structure maps $\eta:A\to \Multiplier(C)$ and $\epsilon:C\to A$. We define $\LAdj(C)\coloneqq \phi^\epsilon : C \to  C^\epsilon$, the local adjunction between the $C^*$-correspondences ${}_A C_C$ and ${}_C C^\epsilon_A$ constructed above.
\end{definition}

The following lemma completes the proof of Theorem \ref{thm:frob-to-ladj}:

\begin{lemma}\label{lem:LAdjFrob}
The canonical isomorphism of $C^*$-algebras $\psi : C\to \Compact_C(C)$ from Theorem \ref{thm:compacts}(2) is an isomorphism of Frobenius $C^*$-algebras $C \xrightarrow{\cong} \Frob\LAdj(C)$.
\end{lemma}

\begin{proof}
Since the $A$-action on the $C^*$-correspondence $C$---and hence, the action of $A$ as multipliers of $\Compact_C(C)$---was defined by means of the homomorphism $\eta$, the map $\psi$ is clearly an $A$-bimodule homomorphism. Thus we are left to prove that $\epsilon_\phi \circ \psi = \epsilon$ as maps $C\to A$, where we have simplified the notation by writing $\phi$ instead of $\phi^\epsilon$. Every element of the $C^*$-algebra $C$ can be written as a product $c_1 c_2$, and on such elements we have
\[
\epsilon_{\phi}\circ \psi(c_1c_2) = \epsilon_{\phi}(\ket{c_1}\bra{c_2^*}) = \langle \phi(c_1)\, |\, \phi(c_2)\rangle_{C^\epsilon} = \epsilon(c_1c_2)
\]
as required.
\end{proof}

\subsection{An endomorphism-ring theorem for Frobenius $C^*$-algebras}

Kasch's endomorphism-ring theorem for Frobenius algebras  asserts that if $C$ is a Frobenius algebra over $A$ (in the algebraic sense of Definition \ref{def:Frobenius-algebra}), then $\End_A(C_A)$---the ring of endomorphisms of $C$, considered as a right $A$-module---is a Frobenius algebra  over $C$ \cite[Satz 1]{Kasch}. There is an analogous result for Frobenius $C^*$-algebras:

\begin{corollary}\label{cor:endomorphism}
Let $C$ be a Frobenius $C^*$-algebra over $A$,  and let $C^\epsilon$ be the $C^*$-correspondence from $C$ to $A$ constructed in Definition \ref{def:C-epsilon}. Then the $C^*$-algebra $\Compact_A(C^\epsilon)$ is, in a canonical way, a Frobenius $C^*$-algebra over $C$.
\end{corollary}

\begin{proof}
The map $(\phi^{\epsilon})^{-1}:C^\epsilon\to C$ (cf.~Definition \ref{def:LAdj} and the paragraph preceding it) defines a local adjunction between the $C^*$-correspondences ${}_C C^\epsilon_A$ and ${}_A C_C$. Applying Proposition \ref{prop:ladj-to-frob} to this local adjunction yields a Frobenius $C^*$-algebra structure, over $C$, on the $C^*$-algebra $\Compact_A(C^\epsilon)$.
\end{proof}

\section{Isomorphism of Frobenius $C^*$-algebras and of local adjunctions}\label{sec:isos}

We have seen that Frobenius $C^*$-algebras can be constructed from local adjunctions of $C^*$-correspondences, and vice versa. In this section we shall sharpen these constructions into a one-to-one correspondence between isomorphism classes of Frobenius $C^*$-algebras on the one hand, and of local adjunctions on the other hand. Two local adjunctions will be considered isomorphic, loosely speaking, if they differ by an equivalence. The formal definition will require a little preparation.

Firstly, if $F$ is a Hilbert $B$-module then we define
\[
I_F\coloneqq \overline{\lspan}\{ \langle f_1\, |\, f_2\rangle\in B\ |\ f_1,f_2\in F\},
\]
a two-sided ideal in $B$. Now let $B$ and $B'$ be $C^*$-algebras, and let $I\subseteq B$ and $I'\subseteq B'$ be ideals. A $C^*$-correspondence ${}_{B'} G_{B}$ from $B'$ to $B$ is a \emph{Morita equivalence from $I'$ to $I$} if $I_G=I$, and if  the left action of $B'$ on $G$ restricts to an isomorphism of $C^*$-algebras $I'\xrightarrow{\cong}\Compact_{B}(G)$.  If $G$ is such an equivalence, then $G^\star\coloneqq \Compact_{B}(G,B)$ is a Morita equivalence from $I$ to $I'$: the $B$-$B'$-bimodule structure on $G^\star$ is given by $b\bra{g}b'\coloneqq \bra{{b'}^* g {b}^*}$, and the $B'$-valued inner product is $\big\langle \bra{g_1}\, \big|\, \bra{g_2}\big\rangle \coloneqq \ket{g_1}\bra{g_2}\in \Compact_{B}(G)\cong I'\subseteq B'$. The map ${}_{B'}G_B\to {}_B G^\star_{B'}$ sending $g\mapsto \bra{g}$ is a local adjunction; and from each local adjunction $\phi:{}_A F'_{B'}\to {}_{B'} E'_A$ we obtain a new local adjunction ${}_A(F'\otimes^{\h}_{B'} G)_{B} \to {}_{B}(G^\star\otimes^{\h}_{B'} E')_A$ via the formula $f\otimes g\mapsto \bra{g}\otimes\phi(f)$. As a final piece of terminology, let us say `a local adjunction \emph{over $A$}' to mean a local adjunction of $C^*$-correspondences of the form $\phi: {}_A F_B\to {}_B E_A$, for some second $C^*$-algebra $B$.

\begin{definition}
Let $\phi:{}_A F_B \to {}_B E_A$ and $\phi' : {}_A F'_{B'} \to {}_{B'} E'_A$ be local adjunctions over $A$. An \emph{isomorphism} $\phi' \xrightarrow{\cong} \phi$ is a triple $(G,\alpha,\beta)$, where $G$ is a $C^*$-correspondence from $B'$ to $B$ that is a Morita equivalence from $I_{F'}$ to $I_{F}$, and where $\alpha : F'\otimes^{\h}_{B'} G\xrightarrow{\cong} F$ and $\beta : G^\star\otimes^{\h}_{B'} E'\to E$ are unitary isomorphisms of $C^*$-correspondences satisfying
\[
\phi\alpha( f'\otimes g  ) = \beta\left( \bra{g}\otimes\phi'(f')\right)
\]
for all $f'\in F'$ and $g\in G$.
\end{definition}

\begin{theorem}\label{thm:iso}
Let $A$ be a $C^*$-algebra.
\begin{enumerate}[\rm(1)]
\item If $C$ is a Frobenius $C^*$-algebra over $A$ then $C \cong \Frob\LAdj(C)$.
\item If $\phi$ is a local adjunction over $A$  then $\phi \cong \LAdj\Frob(\phi)$.
\item If $C$ and $D$ are Frobenius $C^*$-algebras over $A$, then we have $C \cong D$ as Frobenius $C^*$-algebras if and only if $\LAdj(C)\cong \LAdj(D)$ as local adjunctions.
\end{enumerate}
\end{theorem}

\begin{proof}
Part (1) was proved as part of Theorem \ref{thm:frob-to-ladj}.

For part (2), let $\phi: {}_A F_B \to {}_B E_A$ be a local adjunction. We have $\LAdj\Frob(\phi)=\phi^\epsilon : C\to C^\epsilon$,  $c\mapsto c^*$, where $C=\Compact_B(F)$ is considered as a $C^*$-correspondence from $A$ to $C$ in the canonical way; and where $C^\epsilon = C$ considered as a $C^*$-correspondence from $C$ to $A$ with the obvious bimodule structure and with the inner product $\langle c_1\, |\, c_2\rangle^{\epsilon} \coloneqq \epsilon(c_1^*c_2)$, where $\epsilon(\ket{f_1}\bra{f_2})=\langle \phi(f_1)\, |\, \phi(f_2)\rangle$.

Now let $G=F$, considered as a $C^*$-correspondence from $C$ to $B$: the right Hilbert $B$-module structure on $G$ is the given one on $F$, and the left $C$-module structure on $F$ is the canonical action of $C=\Compact_B(F)$. This $C^*$-correspondence ${}_C G_B$ is tautologically a Morita equivalence from $C$ to $I_F$. 

Let $\alpha : C\otimes^{\h}_C G \to F$ be the map $\alpha(c\otimes f)\coloneqq c(f)$, and let $\beta: G^\star \otimes^{\h}_C C^\epsilon \to E$ be the map $\beta\left(\bra{f} \otimes c\right)\coloneqq \phi(c^*(f))$. The map $\beta$ is surjective, by  part (4) of Theorem \ref{thm:compacts}; the computation
\[
\begin{aligned}
 & \Big\langle \beta(\bra{f_1}\otimes c_1)\, \Big|\, \beta(\bra{f_2}\otimes c_2)\Big \rangle  = \Big\langle \phi(c_1^*(f_1))\, \Big|\, \phi(c_2^*(f_2))\Big\rangle  = \epsilon\Big( \ket{c_1^*(f_1)}\bra{c_2^*(f_2)}\Big) \\
& = \epsilon\Big( c_1^*  \ket{f_1}\bra{f_2}c_2\Big) = \Big\langle c_1\, \Big| \, \ket{f_1}\bra{f_2} c_2\Big\rangle^{\epsilon} = \Big\langle \bra{f_1}\otimes c_1\, \Big|\, \bra{f_2}\otimes c_2\Big\rangle
\end{aligned}
\]
shows that $\beta$ preserves the $A$-valued inner products; and the intertwining property of $\phi$ shows that $\beta$ is a $B$-$A$-bimodule map. Thus $\beta$ is a unitary isomorphism of correspondences. A similar argument shows that $\alpha$ is a unitary isomorphism. Finally, the computation
\[
\beta\left( \bra{f}\otimes \phi^{\epsilon}(c) \right) = \phi(c(f)) = \phi \alpha (c\otimes f)
\]
shows that the triple $(G,\alpha,\beta)$ determines an isomorphism of local adjunctions $\phi^\epsilon \to \phi$, and this completes the proof of part (2).

For part (3), let $C$ and $D$ be Frobenius $C^*$-algebras over $A$, with structure maps  $\epsilon:C\to A$ and $\theta:D \to A$. First  suppose that $\rho:D\to C$ is an isomorphism of Frobenius $C^*$-algebras. Let ${}_{D} G_{C}$ be the $C^*$-correspondence defined by $G=C$, with the canonical Hilbert $C$-module structure, and with left $D$-action given by the isomorphism of $C^*$-algebras $\rho : D\to C\cong \Compact_{C}(G)$. Then $G$ is a Morita equivalence from $D$ to $C$. Let $\alpha : D\otimes^{\h}_{D} G \to C$ be the map $\alpha(d\otimes c)=\rho(d)c$, and let $\beta: G^\star \otimes^{\h}_{D} {D}^{\theta} \to C^\epsilon$ be the map $\beta(\bra{c}\otimes d) = c^*\rho( d)$. Straightforward computations show that $\alpha$ and $\beta$ are unitary isomorphisms of correspondences satisfying $\alpha(d\otimes c)^* = \beta( \bra{c}\otimes {d}^*)$, and so $(G,\alpha,\beta)$ is an isomorphism of local adjunctions $\LAdj(D)\to \LAdj(C)$. 

On the other hand, suppose that $(G,\alpha,\beta)$ is an isomorphism of local adjunctions $\LAdj(D)\to \LAdj(C)$: thus ${}_{D} G_{C}$ is a Morita equivalence from $D$ to $C$, and $\alpha: {}_A(D\otimes^{\h}_{D} G)_C\to {}_A C_C$ and $\beta: {}_C(G^\star \otimes^{\h}_{D} {D}^{\theta})_A \to {}_C C^\epsilon_A$ are unitary isomorphisms of $C^*$-corresondences satisfying $\alpha(d\otimes g)^* = \beta( \bra{g}\otimes {d}^*)$. Let $\rho:D\to C$ be the isomorphism of $C^*$-algebras defined as the following composition:
\begin{equation}\label{eq:rho}
D \xrightarrow[\cong]{\text{action on $G$}} \Compact_C(G) \xrightarrow[\cong]{k\mapsto \id_{C'}\otimes k} \Compact_C(D\otimes^{\h}_{D} G) \xrightarrow[\cong]{k\mapsto \alpha\circ k \circ \alpha^*} \Compact_C(C)\xrightarrow[\cong]{} C.
\end{equation}

The fact that $\alpha$ is an $A$-module map ensures that $\rho$ is an $A$-bimodule map. To verify the relation $\epsilon\circ\rho=\theta$, let $\zeta:\Compact_C(D\otimes^{\h}_{D} G)\to D$ be the inverse of the isomorphism $D\xrightarrow{\cong} \Compact_C(G)\xrightarrow{\cong}\Compact_C(D\otimes^{\h}_{D} G)$ appearing in the definition \eqref{eq:rho} of $\rho$. Our goal, then, is to show that $\theta\zeta (k) = \epsilon( \alpha k \alpha^*)$ for all $k\in \Compact_{D}(D\otimes^{\h}_{D}G)$. For $d_1,d_2\in D$ and $g_1,g_2\in G$ we have
\[
\theta\zeta\left( \ket{d_1\otimes g_1}\bra{d_2\otimes g_2}\right) = \big\langle \bra{g_1}\otimes {d_1}^*\, \big|\, \bra{g_2}\otimes {d_2}^* \big\rangle,
\]
the inner product on the right being the $A$-valued inner product on $G^\star\otimes^{\h}_{D} {D}^{\theta}$. Using the fact that $\beta:G^\star\otimes^{\h}_{D} {D}^{\theta}\to C^\epsilon$ is a unitary isomorphism satisfying $\beta(\bra{g}\otimes {d}^*) = \alpha(d\otimes g)^*$, we find
\[
\begin{aligned}
 \theta\zeta\left( \ket{d_1\otimes g_1}\bra{d_2\otimes g_2}\right)&  = \big\langle \beta(\bra{g_1}\otimes d_1^*)\, \big|\, \beta(\bra{g_2}\otimes d_2^* )\big\rangle \\ & = \big\langle \alpha(d_1 \otimes g_1)^*\, \big|\, \alpha(d_2 \otimes g_2)^* \big\rangle \\
& = \epsilon\left(  \alpha(d_1 \otimes g_1)\alpha(d_2\otimes g_2)^*\right)  \\
& = \epsilon\left( \alpha \circ \ket{d_1\otimes g_1}\bra{d_2\otimes g_2} \circ \alpha^*\right)
\end{aligned}
\]
as required. Thus $\rho$ is an isomorphism of Frobenius $C^*$-algebras over $A$.
\end{proof}

\bibliographystyle{alpha}
\bibliography{frobenius}

\end{document}